# NORMAL EXTENSIONS OF A SINGULAR MULTİPOİNT DIFFERENTIAL OPERATOR FOR FİRST ORDER

by

**Z. I. ISMAILOV** and **R.ÖZTÜRK MERT**

Karadeniz Technical University, Faculty of Sciences, Department of Mathematics, 61080, Trabzon, TURKEY; e-mail: zameddin@yahoo.com;

**Abstract**

In this work, firstly in the direct sum of Hilbert spaces of vector-functions $L^2(H,(-\infty,a_1))\oplus L^2(H,(a_2,b_2))\oplus L^2(H,(a_3,+\infty))$, $-\infty < a_1 < a_2 < b_2 < a_3 < +\infty$ all normal extensions of the minimal operator generated by linear singular multipoint formally normal differential expression $l = (l_1, l_2, l_3), l_k = \frac{d}{dt} + A_k$ with a selfadjoint operator coefficient $A_k$ $k = 1,2,3$ in any Hilbert space $H$, are described in terms of boundary values. Later structure of the spectrum of these extensions is investigated.

1. **Introduction**

Many problems arising in the modelling of processes of multi-particle quantum mechanics,in the quantum field theory,in the multipoint boundary value problems for the differential equations,in the physics of rigid bodies and ets support to study normal extension of formally narmal differential operators in direct sum of Hilbert spaces([1-3]).

The general theory of normal extensions of formally normal operators in any Hilbert space have been investigated by many mathematicians ([4-9]). Applications of this theory to two point differential operators in Hilbert space of functions are continied today even(10-17]).

İt is clear that, the minimal operators $L_0(1,0,0) = L_{10}\oplus 0\oplus 0$ and $L_0(0,0,1) = 0\oplus 0\oplus L_{30}$ generated by differential expressions in forms $\left(\frac{d}{dt}+A_1,0,0\right)$ and $(0,0,\frac{d}{dt}+A_3)$ in the Hilbert spaces of vector-functions $L_2(1,0,0) = L_2(H,(-\infty,a_1))\oplus 0\oplus 0$, $L_2(0,0,1) = 0\oplus 0\oplus L_2(H,(a_3,\infty))$ respectively,where
$A_1 = A_1^* \leq 0, A_3 = A_3^* \geq 0, -\infty < a < b < \infty$ ,are maximal formally normal,consequently have not normal extensions.But the minimal operator $L_0(0,1,0) = 0\oplus L_{20}\oplus 0$ generated by differential expression in form $(0,\frac{d}{dt}+A_2,0)$ in the Hilbert spaces of vector-functions $L_2(0,1,0) = 0\oplus L_2(H,(a_2,b_2))\oplus 0$ is formally normal,no is maximal.

Unfortunately,in the multipoint situations may occur in different tables in the following sense.Let $B_1$ , $B_2$ and $B_3$ be minimal operators generated by the linear differential expression $\frac{d}{dt}$ in the Hilbert space of functions $L^2(-\infty, a_1)$, $L^2(a_2, b_2)$ and $L^2(a_3, +\infty)$, $-\infty < a_1 < a_2 < b_2 < a_3 < \infty$, respectively. Consequently, $B_1$ and $B_3$ are maximal formally normal operators, but are not a normal extensions. However, direct sum $B_1 \oplus B_2 \oplus B_3$ of operators $B_1$, $B_2$ and $B_3$ in a direct sum $L^2(-\infty, a_1) \oplus L^2(a_2, b_2) \oplus L^2(a_3, +\infty)$ has a normal extensions. For example, in case $H = \mathbb{C}$ it is can be easy shown that a extension of the minimal operator $B_1 \oplus B_2 \oplus B_3$ with the boundary conditions
$$u_3(a_3) = e^{i\varphi} u_1(a_1), \varphi \in [0, 2\pi),$$
$$u_2(b_2) = e^{i\psi} u_2(a_2), \psi \in [0, 2\pi),$$
$$u = (u_1, u_2, u_3), u_1 \in D(B_1^*), u_2 \in D(B_2^*), u_3 \in D(B_3^*)$$
is a normal in $L^2(-\infty, a_1) \oplus L^2(a_2, b_2) \oplus L^2(a_3, +\infty)$.

In general case of $H$ direct sum $L_0(1,1,1) = L_{10} \oplus L_{20} \oplus L_{30}$ of operators $L_{10}$, $L_{20}$ and $L_{30}$ is formally normal, is not maximal, moreover it have normal extensions in the direct sum $L_2(1,1,1) = L_2(H, (-\infty, a_1)) \oplus L^2(H, (a_2, b_2)) \oplus L_2(H, (a_3, \infty))$.

In singular cases although to date no investigation so far has not been. But the physical and technical processes, many of the problems resulting from examination of the solution is of great importance for the singular cases.

In this paper will be considered a linear multipoint differential-operator expression
$$l = (l_1, l_2, l_3), l_k = \frac{d}{dt} + A_k, k = 1,2,3$$
in the direct sum of Hilbert spaces of vector-functions $L_2(1,1,1)$, where $A_1 = A_1^* \leq 0$, $A_2 = A_2^* \geq 0$, $A_3 = A_3^* \geq 0 - \infty < a_1 < a_2 < b_2 < a_3 < \infty$.

In this work in second section by the method of J.W.Calkin- M.L.Gorbachuk theory all normal extensions of the minimal operator generated by singular multipoint formally normal differential expression for first order $l(.)$ in the direct sum of Hilbert space $L_2(1,1,1)$ in terms of boundary values are described. In third section the spectrum of such extensions is researched.

## 2. Description of Normal Extensions

Let $H$ be a separable Hilbert space and $a_1, a_2, b_2, a_3 \in R, a_1 < a_2 < b_2 < a_3$. In the Hilbert space $L_2(1,1,1)$ of vector-functions consider the following linear multipoint differential expression in form
$$l(u) = (l_1(u_1), l_2(u_2), l_3(u_3)) = (u_1' + A_1 u_1, u_2' + A_2 u_2, u_3' + A_3 u_3), u = (u_1, u_2, u_3),$$

where $A_k: D(A_k) \subset H \to H, k = 1,2,3$ is linear selfadjoint operators in $H$ and $A_1 = A_1^* \leq 0, A_2 = A_2^* \geq 0, A_3 = A_3^* \geq 0$. In the linear manifold $D(A_k) \subset H$ introduce the inner product in form
$$(f, g)_{k,+} := (A_k f A_k, g)_H + (f, g)_H, \ f, g \in D(A_k), k = 1,2,3.$$
For the $k = 1,2,3$, $D(A_k)$ is a Hilbert space under the positive norm $\|\cdot\|_{k,+}$ respect to Hilbert space $H$. It is denoted by $H_{k,+}$. Denote the $H_{k,-}$ a Hilbert space with negative norm. It is clear that a operator $A_k$ is continuous from $H_{k,+}$ to $H$ and it's adjoint operator $\tilde{A}_k: H \to H_{k,-}$ is a extension of the operator $A_k$. On the other hand, the operator $\tilde{A}_k: D(\tilde{A}_k) = H \subset H_{k,-1} \to H_{k,-1}$ is a linear selfadjoint.

Now in the direct sum $L_2(1,1,1)$ define by
$$\tilde{l}(u) = (\tilde{l}_1(u_1), \tilde{l}_2(u_2), \tilde{l}_3(u_3)), \tag{2.1}$$
where $u = (u_1, u_2, u_3)$ and $\tilde{l}_1(u_1) = u_1' + \tilde{A}_1 u_1$, $\tilde{l}_2(u_2) = u_2' + \tilde{A}_2 u_2, \tilde{l}_2(u_2) = u_2' + \tilde{A}_2 u_2, \tilde{l}_3(u_3) = u_3' + \tilde{A}_3 u_3$.

The operators defined by $L_0(1,1,1) = L_{10} \oplus L_{20} \oplus L_{30}$ and $L(1,1,1) = L_1 \oplus L_2 \oplus L_3$ in the space $L_2(1,1,1)$ are called minimal(multipoint) and maximal (multipoint) operators generated by the differential expression (2.1), respectively.

Here it is described all normal extensions of the minimal operator $L_0(1,1,1)$ in $L_2(1,1,1)$ in terms of the boundary values.

Note that a space of boundary values has an important role in the theory extensions of the linear symmetric differential operators ( [18] ) which from it will be used in last investigation.

Let $B: D(B) \subset \mathcal{H} \to \mathcal{H}$ be a closed densely defined symmetric operator in the Hilbert space $\mathcal{H}$, having equal finite or infinite deficiency indices. A triplet $(\mathfrak{H}, \gamma_1, \gamma_2)$, where $\mathfrak{H}$ is a Hilbert space, $\gamma_1$ and $\gamma_2$ are linear mappings of $D(B^*)$ into $\mathfrak{H}$, is called a space of boundary values for the operator $B$ if for any $f, g \in D(B^*)$

$$(B^*f, g)_\mathcal{H} - (f, B^*g)_\mathcal{H} = (\gamma_1(f), \gamma_2(g))_\mathfrak{H} - (\gamma_2(f), \gamma_1(g))_\mathfrak{H},$$

while for any $F, G \in \mathfrak{H}$, there exists an element $f \in D(B^*)$, such that $\gamma_1(f) = F$ and $\gamma_2(f) = G$.

Note that any symmetric operator with equal deficiency indices have at least one space of boundary values ([18]).

Now construct a space of boundary values for the minimal operators $M_0(1,0,1)$ and $M_0(0,1,0)$ generated by linear singular differential expressions of first order in form

$$(m_1(u_1), 0, m_3(u_3)) = (-i\frac{du_1}{dt}, 0, -i\frac{du_3}{dt})$$

$$(0, m_2(u_2), 0) = (0, -i\frac{du_2}{dt}, 0)$$

in the direct sum $L_2(1,0,1)$ and $L_2(0,1,0)$, respectively. Note that the minimal operators $M_0(1,0,1)$ and $M_0(0,1,0)$ are closed symmetric operators in $L_2(1,0,1)$ and $L_2(0,1,0)$ with deficiency indices $(dimH, dimH)$.

**Lemma.2.1.** The triplet $(H, Y_1, Y_2)$, where

$$Y_1: D(M_0^*) \to H, Y_1(u) = \frac{1}{i\sqrt{2}}(u_3(a_3) + u_1(a_1)),$$

$$Y_2: D(M_0^*) \to H, Y_2(u) = \frac{1}{\sqrt{2}}(u_3(a_3) - u_1(a_1)), \quad u = (u_1, 0, u_3) \in D(M_0^*)$$

is a space of boundary values of the minimal operator $M_0(1,0,1)$ in $L_2(1,0,1)$.

**Proof.** For arbitrary $u = (u_1, 0, u_3)$ and $v = (v_1, 0, v_3)$ from $D(M_0^*(1,0,1))$ validity the following equality

$$(M_0^*(1,0,1)u, v)_{L_2(1,0,1)} - (u, M_0^*(1,0,1)v)_{L_2(1,0,1)} = (Y_1(u), Y_2(v))_H - (Y_2(u), Y_1(v))_H$$

can be easily verified. Now give any elements $f, g \in H$. Find the function $u = (u_1, 0, u_3) \in D(M_0^*(1,0,1))$ such that

$Y_1(u) = \frac{1}{i\sqrt{2}}(u_3(a_3) + u_1(a_1)) = f$ and $Y_2(u) = \frac{1}{\sqrt{2}}(u_3(a_3) - u_1(a_1)) = g$,

that is,

$u_1(a_1) = (if - g)/\sqrt{2}$ and $u_3(a_3) = (if + g)/\sqrt{2}$.

If choose these functions $u_1(t), u_2(t)$ in following form

$u_1(t) = \int_{-\infty}^{t} e^{s-a} ds (if - g)/\sqrt{2}, \quad t < a_1, \quad u_3(t) = \int_{t}^{\infty} e^{a_3-t} ds (if + g)/\sqrt{2}, \quad t > a_3,$

then it is clear that $(u_1, u_2) \in D(M_0^*)$ and $Y_1(u) = f$, $Y_2(u) = g$.

**Lemma 2.2**. The triplet $(H, \Gamma_1, \Gamma_2)$,
$$\Gamma_1: D(M_0^*(0,1,0)) \to H, \ \Gamma_2(u) = \frac{1}{i\sqrt{2}}(u_2(b_2) + u_2(a_2)),$$
$$\Gamma_2: D(M_0^*(0,1,0)) \to H, \Gamma_2(u) = \frac{1}{\sqrt{2}}(u_2(b_2) - u_2(a_2)),$$
$$u = (0, u_2, 0) \in D(M_0^*(0,1,0))$$
is a space of boundary values of the minimal operator $M_0(0,1,0)$ in the direct sum $L_2(0,1,0)$.

**Theorem 2.3.** If the minimal operators $L_{10}, L_{20}$ and $L_{30}$ are formally normal, then are true
$$D(L_{10}) \subset W_2^1(H, (-\infty, a_1)), A_1 D(L_{10}) \subset L_2(H, (-\infty, a_1)),$$
$$D(L_{20}) \subset W_2^1(H, (a_2, b_2)), A_2 D(L_{20}) \subset L_2(H, (a_2, b_2)),$$
$$D(L_{30}) \subset W_2^1(H, (a_3, \infty)), A_3 D(L_{30}) \subset L_2(H, (a_3, \infty))$$

**Proof.** Indeed, in this case for each $u_1 \in D(L_{10}) \subset D(L_{10}^*)$ is true
$u_1' + A_1 u_1 \in L_2(H, (-\infty, a))$ and $u_1' - A_1 u_1 \in L_2(H, (-\infty, a))$, hence
$u_1' \in L_2(H, (-\infty, a))$ and $A_1 u_1 \in L_2(H, (-\infty, a))$, i.e.,
$$D(L_{10}) \subset W_2^1(H, (-\infty, a)) \quad and \quad A_1 D(L_{10}) \subset L_2(H, (-\infty, a)).$$
In similar way it is proved the second and third parts of theorem.
On the other hand can be easily established the following result.

**Lemma 2.4.** Every normal extension of $L_0(1,1,1)$ in $L_2(1,1,1)$ is a direct sum of normal extensions of the minimal operator $L_0(1,0,1) = L_{10} \oplus 0 \oplus L_{30}$ in $L_2(1,0,1) = L_2(H, (-\infty, a_1)) \oplus 0 \oplus L_2(H, (a_3, \infty))$ and minimal operator $L_0(0,1,0) = 0 \oplus L_{20} \oplus 0$ in $L_2(0,1,0) = 0 \oplus L_2(H, (a_2, b_2)) \oplus 0$.

Finally, using the method in [10,18] and Lemmas 2.1 and 2.2 can be established the following result.

**Theorem 2.5.** Let $(-A_1)^{1/2} W_2^1(H, (-\infty; a_1)) \subset W_2^1(H, (-\infty, a_1))$,
$$A_2^{1/2} W_2^1(H, (a_2, b_2)) \subset W_2^1(H, (a_2, b_2)),$$
and $\quad A_3^{1/2} W_2^1(H, (a_3, \infty)) \subset W_2^1(H, (a_3, \infty))$.

Each normal extention $\tilde{L}$ of the minimal operator $L_0$ in the Hilbert space $L_2(1,1,1)$ is generated by differential expression (2.1) and boundary contitions
$$u_3(a_3) = W_1 u_1(a_1), u_1(a_1) \in \ker(-A_1)^{1/2}, u_3(a_3) \in \ker A_3^{1/2}, \quad (2.3)$$
$$u_2(b_2) = W_2 u_1(a_2), \quad (2.4)$$
where $W_1, W_2: H \to H$ is a unitary operators. Moreover, the unitary operators $W_1, W_2$ in $H$ are determined by the extension $\tilde{L}$, i.e. $\tilde{L} = L_{W_1 W_2}$ and vice versa.

**Corollary 2.6.** If at least one of the operators $A_1$ and $A_3$ is one-to-one mapping in $H$, then minimal operator $L_0(1,1,1)$ is maximally formal normal in $L_2(1,1,1)$.

**Corollary 2.7.** If there exists at least one normal extension of the minimal operator $L_0(1,1,1)$, then
$$dimker(-A_1)^{1/2} = dimker A_3^{1/2} > 0.$$

### 3. The Spectrum of the Normal Extensions

In this section the structure of the spectrum of the normal extension $L_{W_1 W_2}$ in $L_2(1,1,1)$ will be investigated. In this case by the Lemma 2.4. it is clear that
$$L_{W_1 W_2} = L_{W_1} \oplus L_{W_2},$$

where $L_{W_1}$ and $L_{W_2}$ are normal extensions of the minimal operators $L_0(1,0,1)$ and $L_0(0,1,0)$ in the Hilbert spaces $L_2(1,0,1)$ and $L_2(0,1,0)$ respectively.

Later will be assumed that $A_1 = A_1^* \leq 0, A_2 = A_2^* \geq 0, A_3 = A_3^* \geq 0$ and $0 \in \sigma_p((-A_1)^{1/2}) \cap \sigma_p(A_3^{1/2})$.

First of all, we have to prove the following result.

**Theorem 3.1:** The point spectrum of any normal extension $L_{W_1}$ of the minimal operator $L_0(1,0,1)$ in the Hilbert space $L_2(1,0,1)$ is empty, i.e.,
$$\sigma_p(L_{W_1}) = \emptyset.$$

**Proof:** In first consider the following problem for the spectrum of the normal extension $L_{W_1}$ of the minimal operator $L_0(1,0,1)$ in the Hilbert space $L_2(1,0,1)$
$$L_{W_1} u = \lambda u, \lambda = \lambda_r + i\lambda_i \in \mathbb{C}, u = (u_1, 0, u_3) \in L^2(1,0,1),$$
that is,
$$\tilde{l}_1(u_1) = u_1' + \tilde{A}_1 u_1 = \lambda u_1, \quad u_1 \in L^2(H, (-\infty, a_1)),$$
$$\tilde{l}_3(u_3) = u_3' + \tilde{A}_3 u_3 = \lambda u_3, \quad u_3 \in L^2(H, (a_3, +\infty)), \quad \lambda \in R$$
$$u_3(a_3) = W_1 u_1(a_1), u_1(a_1) \in \ker(-A_1)^{1/2}, u_3(a_3) \in \ker A_3^{1/2}.$$

The general solution of this problem is

$$u_1(\lambda; t) = e^{-(\tilde{A}_1 - \lambda)(t-a_1)} f_1^*, \quad t < a_1, f_1^* \in H_{-1/2}(-A_1),$$
$$u_3(\lambda; t) = e^{-(\tilde{A}_3 - \lambda)(t-a_3)} f_3^*, \quad t > a_3, f_3^* \in H_{-1/2}(A_3)$$
$$f_3^* = W_1 f_1^*, \quad f_1^*, f_3^* \in H, f_1^* = u_1(\lambda; a_1), f_3^* = u_3(\lambda; a_3)$$

Since $0 \in \sigma_p((-A_1)^{1/2}) \cap \sigma_p(A_2^{1/2})$ and $(-A_1)^{1/2} f_1^* = 0$, $A_2^{1/2} f_3^* = 0$, then we have
$$u_1(\lambda; t) = e^{\lambda(t-a)} f_1^*, t < a, f_1^* \in H_{-1/2}((-A_1))$$
$$u_2(\lambda; t) = e^{\lambda(t-b)} f_3^*, t > b, f_3^* \in H_{-1/2}(A_3),$$
$$f_3^* = W_1 f_1^*, \; , f_1^* = u_1(\lambda; a_1), f_3^* = u_3(\lambda; a_3).$$

In order to $u_1(\lambda; .) \in L_2(H, (-\infty, a_1))$ and $u_2(\lambda; .) \in L_2(H, (a_3, \infty))$ necessary and sufficient condition is $\lambda_r \geq 0$ and $\lambda_r \leq 0$ respectively. Hence $\lambda_r = 0$.

Consequently,
$$u_1(\lambda; t) = e^{i\lambda_i(t-a_1)} f_1^*, t < a_1,$$
$$u_2(\lambda; t) = e^{i\lambda_i(t-a_3)} f_3^*, t > a_3, f_3^* = W_1 f_1^*.$$

In this case it is clear that for the $u_1(\lambda; .) \in L_2(H, (-\infty, a_1))$ and $u_2(\lambda; .) \in L_2(H, (a_3, \infty))$ necessary and sufficient conditions are $f_1^* = 0, f_3^* = 0$. From this implies that $u_1 = 0$ and $u_2 = 0$ in $L^2$. Therefore $\sigma_p(L_{W_1}) = \emptyset$.

Since residual spectrum of any normal operators in any Hilbert space is empty, then furthermore the continuous spectrum of the normal extensions $L_{W_1}$ of the minimal operator $L_0(1,0,1)$ in the Hilbert space $L_2(1,0,1)$ is investigated.

**Theorem 3.2:** Continuous spectrum of the any normal extension $L_{W_1}$ of the minimal operator $L_0(1,0,1)$ in the Hilbert space $L_2(1,0,1)$ is $\sigma_c(L_{W_1}) = i\mathbb{R}$.

**Proof**. Assume that $\lambda \in \sigma_c(L_{W_1})$. Then on the important theorem for the spectrum of normal operators [19], that is,
$$\sigma(L_{W_1}) \subset \sigma(ReL_{W_1}) + i\sigma(İmL_{W_1}),$$
it is obtained that
$$\lambda_r \in \sigma(ReL_{W_1}), \lambda_i \in \sigma(İmL_{W_1}).$$
From this imply that $\lambda_r \in \sigma(A_1), \lambda_r \in \sigma(A_3)$, hence on the conditions to the operators $A_1$ and $A_3$ we have $\lambda_r = 0$. On the other hand from the proof of previes theorem it is to see that $\ker(L_{W_1} - \lambda) = \{0\}$ for any $\lambda \in \mathbb{C}$. Consequently, $\sigma_c(L_{W_1}) \subset i\mathbb{R}$.

Furthermore, it is clear that for the $\lambda = i\lambda_i \in \mathbb{C}$ the general solution of the boundary value problem

$$u_1' + A_1 u_1 = i\lambda_i u_1 + f_1, \quad u_1, f_1 \in L_2(H, (-\infty, a_1)),$$
$$u_3' + A_3 u_3 = i\lambda_i u_3 + f_3, \quad u_3, f_3 \in L_2(H, (a_3, \infty)), \lambda_i \in R,$$
$$u_3(a_3) = W_1 u_1(a_1), u_1(a_1) \in \ker(-A_1)^{1/2}, u_3(a_3) \in \ker A_3^{1/2}$$

İs in form
$$u_1(i\lambda_i; t) = e^{-(A_1 - i\lambda_i)(t - a_1)} f_{i\lambda_i} - \int_t^{a_1} e^{-(A_1 - i\lambda_i)(t - s)} f_1(s) ds, \quad t < a_1,$$

$$u_3(i\lambda_i; t) = e^{-(A_3 - i\lambda_i)(t - a_3)} g_{i\lambda_i} + \int_{a_3} e^{-(A_3 - i\lambda_i)(t - s)} f_3(s) ds, \quad t > a_3,$$

$$g_{i\lambda_i} = W_1 f_{i\lambda_i}.$$

In this case it is true
$e^{-(A_1 - i\lambda_i)(t - a_1)} f_{i\lambda_i} \in L_2(H, (-\infty, a_1))$, $e^{-(A_2 - i\lambda_i)(t - a_3)} g_{i\lambda_i} \in L_2(H, (a_3, \infty))$ for any $g_{i\lambda_i}, f_{i\lambda_i} \in H$.

Here if choose $f_1(t) = e^{i\lambda_i t} e^{-(t - a_1)} f^*, f^* \in \ker(-A_1)^{1/2}, t < a_1$, then
$$\int_t^{a_1} e^{-(A_1 - i\lambda_i)(t - s)} f_1(s) ds = e^{-i\lambda_i t} \int_t^{a_1} e^{-(s - a_1)} f^* ds = e^{-i\lambda_i t}(e^{-(t - a_1)} - 1) f^*, t < a_1.$$

Therefore
$$\int_{-\infty}^{a_1} \| e^{-i\lambda_i t}(e^{-(t - a_1)} - 1) f^* \|^2 dt = \int_{-\infty}^{a_1} \| e^{-i\lambda_i t}(e^{-(t - a_1)} - 1) f^* \|^2 dt$$
$$= \int_{-\infty}^{a_1} (e^{-2(t - a_1)} - 2e^{-(t - a_1)} + 1) dt \| f^* \|^2 = \infty$$

Consequently, for the such $f_1(t) \in L_2(H, (-\infty, a_1))$, $u_1(i\lambda_i; t) \notin L_2(H, (-\infty, a_1))$. This means that for any $\lambda \in \mathbb{C}$, a operator $L_{W_1} - \lambda$ is one-to-one in $L_2(1,0,1)$, but it is not onto transformation. On the other hand, since resudial spectrum $\sigma_r(L_{W_1})$ is empty, then it is implies that $\sigma(L_{W_1}) = \sigma_c(L_{W_1}) = i\mathbb{R}$.

Furthermore, here spectrum of normal extensions $L_{W_2}$ of the minimal operator $L_0(0,1,0)$ in $L_2(1,0,1)$ will be investigated.

**Theorem 3.3.** The spectrum of the normalextension $L_{W_2}$ of the minimal operator $L_0(0,1,0)$ in the Hilbert space $L_2(0,1,0)$ is in form

$$\sigma(L_{W_2}) = \Big\{\lambda \in \mathbb{C} : \lambda = \frac{1}{a_2 - b_2}(\ln|\mu| + i \arg\mu + 2n\pi i), n \in \mathbb{Z}, \mu \in \sigma(W_2^* e^{-\tilde{A}_2(b_2 - a_2)}), 0 \le \arg\mu < 2\pi\Big\}$$

**Proof.** The general solution of the following problem to spectrum of the normal extension $L_{W_2}$,

$$\tilde{l}_2(u_2) = u_2' + \tilde{A}_2 u_2 = \lambda u_2 + f_2, \quad u_2, f_2 \in L^2(H,(a_2,b_2))$$
$$u_2(b_2) = W_2 u_2(a_2), \quad \lambda \in \mathbb{C}$$

is in form

$$u_2(t) = e^{-(\tilde{A}_2-\lambda)(t-a_2)} f_2^* + \int_{a_2}^{t} e^{-(\tilde{A}_2-\lambda)(t-s)} f_2(s)ds, a_2 < t < b_2, f_2^* \in H_{-1/2}(A_2)$$

$$(e^{-\lambda(b_2-a_2)} - W_2^* e^{-\tilde{A}_2(b_2-a_2)}) f_2^* = W_2^* e^{-\lambda(b_2-a_2)} \int_{a_2}^{b_2} e^{-(\tilde{A}_2-\lambda)(b_2-s)} f_2(s)ds$$

From this it is implies that $\lambda \in \sigma(L_{W_2})$ if and only if $\lambda$ is a solution of the equation $e^{-\lambda(b_2-a_2)} = \mu$, where $\mu \in \sigma(W_2^* e^{-\tilde{A}_2(b_2-a_2)})$. From this it is obtained that
$\lambda = \frac{1}{a_2-b_2}(\ln|\mu| + i\,arg\mu + 2n\pi i), n \in \mathbb{Z}, \mu \in \sigma(W_2^* e^{-\tilde{A}_2(b_2-a_2)})$.

**Theorem 3.4.** For the spectrum $\sigma(L_{W_1 W_2})$ of any normal extension $L_{W_1 W_2} = L_{W_1} \oplus L_{W_2}$ is true
$$\sigma_p(L_{W_1 W_2}) = \sigma_p(L_{W_2}), \sigma_c(L_{W_1 W_2}) = \{[\sigma_p(L_{W_2})]^c \cap [i\mathbb{R}]\} \cup \sigma_c(L_{W_2})$$

**Proof.** Validity of this assertion is a simple result of the following claim that a proof of which it is clear. If $S_1$ and $S_2$ two linear closed operators in any Hilbert spaces $H_1$ and $H_2$ respectively, then

$$\sigma_p(S_1 \oplus S_2) = \sigma_p(S_1) \cup \sigma_p(S_2),$$
$$\sigma_c(S_1 \oplus S_2) = \left(\sigma_p(S_1) \cup \sigma_p(S_2)\right)^c \cap \left(\sigma_r(S_1) \cup \sigma_r(S_2)\right)^c \cap \left(\sigma_c(S_1) \cup \sigma_c(S_2)\right)$$

Later on, note that for the singular differential operators for n-th order in scalar case in the finite interval has been reseached in [20].

**Example 3.5:** Consider the following boundary value problem for the differential operator $L_{\varphi\psi}$

$$L_{\varphi\psi}: \quad \frac{\partial u(t,x)}{\partial t} + sgnt\frac{\partial^2 u(t,x)}{\partial x^2} = f(t,x), |t| > 1, x \in [0,1],$$
$$\frac{\partial u(t,x)}{\partial t} - \frac{\partial^2 u(t,x)}{\partial x^2} + u(t,x) = f(t,x), \quad |t| < 1/2, x \in [0,1],$$

$$u(1,x) = e^{i\varphi}u(-1,x), \varphi \in [0,2\pi),$$
$$u(1/2,x) = e^{i\psi}u(-1/2,x), \psi \in [0,2\pi),$$
$$u_x(t,0) = u_x(t,1) = 0, |t| > 1, |t| < 1/2$$

in the space $L^2((-\infty,-1) \times (0,1)) \oplus L^2((-1/2,1/2) \times (0,1)) \oplus L^2((1,\infty) \times (0,1))$.
In this case it is clear that in the space $L_2(0,1)$ for the operators

$$A_1 = \frac{\partial^2 u(.,x)}{\partial x^2}, x \in [0,1], u_x(.,0) = u_x(.,1) = 0,$$

$$A_2 = -\frac{\partial^2 u(.,x)}{\partial x^2} + u(.,x) \quad, x \in [0,1]\,, u_x(.,0) = u_x(.,1) = 0,$$

$$A_3 = -\frac{\partial^2 u(.,x)}{\partial x^2} \quad, x \in [0,1]\,, u_x(.,0) = u_x(.,1) = 0$$

are true

$$A_1 = A_1^* \leq 0, A_2 = A_2^* \geq 1, A_3 = A_3^* \geq 0 \,, \ker(-A_1)^{1/2} \neq \{0\}, \ker A_3^{1/2} \neq \{0\},$$
$$0 \in \sigma_p((-A_1)^{1/2}) \cap \sigma_p(A_3^{1/2}).$$

On the other hand, since $A_2^{-1} \in \sigma_\infty(L_2(0,1))$, $\sigma(L_\psi) = \sigma_p(L_\psi)$, $\sigma_c(L_\psi) = \emptyset$ and

$$\sigma(L_\psi) = \{\lambda \in \mathbb{C}: \lambda = \ln|\mu| + i\,arg\mu + 2n\pi i, n \in \mathbb{Z}\,, \mu \in \sigma(e^{i\psi}e^{-\tilde{A}_2(b_2-a_2)}), 0 \leq arg\,\mu < 2\pi\} \subset \{\lambda \in \mathbb{C}: \text{Re}\lambda \geq 1\},$$

then

$$[\sigma_p(L_\psi)]^c \cap [i\mathbb{R}] = i\mathbb{R}.$$

Therefore by the Theorem 3.4 it is obtained

$$\sigma_p(L_{\varphi\psi}) = \sigma_p(L_\psi)\,, \sigma_c(L_{\varphi\psi}) = i\mathbb{R}.$$

**Acknowledgment**

We thank Prof.M.L.Gorbachuk(Institute of Mathematics NASU,Kiev,Ukraine) for his advice and enthusiastic support.